\input amstex
\documentstyle{amsppt}
\magnification=\magstep1

\pageheight{9.0truein}
\pagewidth{6.5truein}

\NoBlackBoxes
\TagsAsMath
\TagsOnRight

%%\long\def\ignore#1{INSERT XY-GRAPHIC}
\long\def\ignore#1{#1}

\ignore{
\input xy
\xyoption{matrix}\xyoption{arrow} %% \xyoption{curve}\xyoption{frame}
\def\edge{\ar@{-}}
}

\def\Rmod{R{}\operatorname{-Mod}}

\def\Zmod{{\Bbb Z}{}\operatorname{-Mod}}

\def\NN{{\Bbb N}}
\def\ZZ{{\Bbb Z}}

\def\AmiProc{{\bf 1}}
\def\Auslarge{{\bf 2}}
\def\Chase{{\bf 3}}
\def\CoKa{{\bf 4}}
\def\Duf{{\bf 5}}
\def\GruJen{{\bf 6}}
\def\Herz{{\bf 7}}
\def\Huibiel{{\bf 8}}
\def\HuiOko{{\bf 9}}
\def\HuiZimI{{\bf 10}}
\def\HuiZimII{{\bf 11}}
\def\Koe{{\bf 12}}
\def\Rowen{{\bf 13}}
\def\Schmid{{\bf 14}}
\def\Simpad{{\bf 15}}
\def\Zim{{\bf 16}}

\document
 
\topmatter

\title  Direct products of modules and the pure semisimplicity conjecture. Part II 
\endtitle

\rightheadtext{ Direct products }

\author Birge Huisgen-Zimmermann and Manuel
Saor\'\i n\endauthor

\address Department of Mathematics, University of California, Santa Barbara, CA
93106, USA\endaddress

\email birge\@math.ucsb.edu\endemail

\address Departamento de Matem\'aticas, Universidad de Murcia, 30100
Espinardo-MU, Spain \endaddress

\email msaorinc\@um.es \endemail

\abstract  We prove that the module categories of Noether algebras (i.e.,
algebras module finite over a noetherian center) and affine noetherian PI algebras over
a field enjoy the following product property:  Whenever a direct product $\prod_{n
\in \NN} M_n$ of finitely generated indecomposable modules $M_n$ is a direct sum of
finitely generated objects, there are repeats among the isomorphism types of the
$M_n$.  The rings with this property satisfy the pure semisimplicity conjecture
which stipulates that vanishing one-sided pure global dimension entails finite
representation type. \endabstract     

\endtopmatter

This is a follow-up to a joint article of the firstnamed author with F. Okoh
\cite{\HuiOko}, improving significantly on the main result of that note.  In
rough terms, our aim is to identify the noetherian rings whose module categories display
the same product-decomposition properties as have long been observed in those of finite
dimensional algebras.  So far, such rings are known to include Artin algebras (this is
essentially due to Auslander \cite{\Auslarge, Corollary 3.2}) and commutative noetherian
domains of Krull dimension 1 (as was proved in \cite{\HuiOko}).  As we show here, not
only can all requirements on the Krull dimension be dropped, but rudimentary
commutativity conditions already guarantee the desired product-behavior.

The principal motivation for our interest in decomposition properties of direct
products lies in the fact that they massively impinge on global decomposition
patterns within the pertinent module categories (see, e.g.,
\cite{\GruJen}, \cite{\Zim}, \cite{\HuiZimI}).  So, in particular, it is the lack
of understanding of direct products which is responsible for the fact
that it is still unresolved whether the rings whose left modules split into
finitely generated submodules (the {\it left pure semisimple rings}) are
necessarily of finite representation type.  For a synopsis of the extensive history
of this pursuit, going back to work of Koethe in the 1930's and of Cohen-Kaplansky
in the early 1950's (\cite{\Koe}, \cite{\CoKa}), the reader is referred to
\cite{\HuiZimII},
\cite{\Simpad}, or \cite{\Huibiel}. The following question crystallizes the
remaining difficulties on the road towards a full resolution of the `pure
semisimplicity conjecture' (which, at this point, is believed to fail in general).  
\smallskip  

{\bf Central Problem:}  For which rings $R$ (associative and with identity) does the
following hold?  If $(M_n)_{n \in \NN}$ is any sequence of pairwise non-isomorphic
finitely generated, indecomposable left
$R$-modules, the direct product $\prod_{n \in \NN} M_n$ is not a direct sum of
finitely generated submodules.
\smallskip

In this form, the question was raised by Okoh.  It connects with the
pure semisimplicity problem as follows:  Whenever a class of rings is known to satisfy
the specified product condition, all left pure semisimple members of that class have
finite representation type.  In our opinion, this problem also holds
independent interest, as a yardstick measuring the current level of understanding
of submodule lattices of direct products.  

In \cite{\HuiOko} it was shown that
commutative noetherian domains of Krull dimension 1 satisfy the above product condition.
This led to the conjecture that the same holds true for arbitrary noetherian
PI-rings, which is in keeping with the fact that the left
pure semisimple artinian PI-rings have finite representation type
(see \cite{\Herz} and \cite{\Schmid}).   Here we come close to confirming it.  Namely,
we establish the product condition for arbitrary Noether algebras (rings which are
module finite over a noetherian center), as well as for affine noetherian PI-rings
(i.e., affine noetherian PI algebras over $\ZZ$). 

For the slightly more general statement of our theorem, we recall that a commutative
ring is a Jacobson ring if each prime ideal is
an intersection of maximal ideals.  So, in particular, affine noetherian PI-algebras
over fields or over noetherian domains of Krull dimension 1 fall into the second class
of rings addressed by the following theorem.  

\proclaim{Theorem}  Suppose that $R$ is 

$\bullet$ either a Noether algebra, or else

$\bullet$ an affine noetherian PI-algebra over a noetherian Jacobson ring.

Then $R$ satisfies the product condition of our `central problem', i.e., given any
sequence $(M_n)_{n \in \NN}$ of pairwise non-isomorphic finitely generated,
indecomposable left
$R$-modules, the direct product $\prod_{n \in \NN} M_n$ fails to be a direct
sum of finitely generated components. \endproclaim

\demo{Proof}   Our strategy consists of playing the problem back to Artin algebras,
where it is already resolved.  To that end, we verify in an initial step that both
classes of rings in our claim enjoy the following property:  For any sequence $P_1,
\dots, P_n$ of (not necessarily distinct) left primitive ideals of $R$, the factor ring
$R / (P_1\cdots P_n)$ is an Artin algebra.

First suppose that $R$ is a Noether algebra and so, in particular, a PI-ring. 
Denote the center of $R$ by $C$.  If $P \subseteq R$ is a left primitive ideal,
Kaplansky's Theorem guarantees $P \cap C$ to be a maximal ideal of $C$;
indeed,
$R/P$ is a finite dimensional algebra over a central subfield $K$ containing $C/ (P \cap
C)$, and since the embedding $C/ (P\cap C) \hookrightarrow K$ is integral, $C/ (P\cap
C)$ is a field as well.  Hence, given any finite sequence $P_1, \dots, P_n$ of left
primitive ideals of $R$, the factor ring $C / \bigl( (P_1 \cap C) \cdots (P_n \cap
C) \bigr)$ has Krull dimension zero and is therefore artinian.  So is, a fortiori, the
central subring
$C / \bigl(  (P_1 \cdots P_n) \cap C \bigr)$  of $R/ (P_1
\cdots P_n)$.  

Next suppose that $C$ is a central noetherian Jacobson subring of $R$.  We will use the
following upgrade of Kaplansky's Theorem, which is due to Amitsur and Procesi
(\cite{\AmiProc}  --  see also
\cite{\Rowen} or
\cite{\Duf} for a slick proof due to Duflo):  If
$R$ is an affine PI-algebra over $C$ and $P$ any left primitive ideal of $R$, then
$P\cap C$ is a maximal ideal of $C$, and
$R/P$ is a simple finite dimensional $\bigl( C/ (P \cap C) \bigr)$-algebra.  Now
suppose that $R$ is noetherian in addition, and let $P_1, \dots, P_n$ be left primitive
ideals of $R$. Since the successive factors of the chain $R\supseteq P_n \supseteq
P_{n-1}P_n
\supseteq
\cdots \supseteq (P_1\cdots P_n)$ are finitely generated left modules over the rings
$R/P_i$ and these rings, in turn, are module-finite over $C$, we see that
$R/ (P_1
\cdots P_n)$ is module-finite over
$C/\bigl(  (P_1 \cdots P_n) \cap C \bigr)$, the latter ring being a homomorphic image
of $C/\bigl( (P_1\cap C)\cdots (P_n \cap C) \bigr)$ and hence artinian.

Consequently, the theorem will follow from the ensuing lemma. \qed \enddemo

\proclaim{Lemma}  Suppose that $R$ is a twosided noetherian ring with the property
that, for any finite sequence
$P_1, \dots, P_n$ of (not necessarily distinct) left primitive ideals, the factor ring
$R/ (P_1
\cdots P_n)$ is an Artin algebra.  Then $R$ satisfies the conclusion of the theorem.
\endproclaim

\demo{Proof}  Let $(M_n)_{n \in \NN}$ be a sequence of finitely generated
indecomposable left $R$-modules with $M_i \not\cong M_j$ whenever $i \ne j$.  Start
by observing that, for any nonzero finitely generated left $R$-module $X$, there exists
a left primitive ideal $P \subseteq R$ with $PX \subsetneqq X$.  Indeed, if $Y$ is a
maximal submodule of $X$, the conductor ideal $P = [Y:X]$ satisfies these
requirements.  

We now apply this observation to finitely generated left $R$-modules of the form $A
M_n$, where
$A$ is a twosided ideal.  Our goal is to construct a descending chain $(A_n)_{n \in
\NN}$ of ideals, each term of which is a finite product of left primitive ideals, with
the following property:  For $n \in \NN$, either $A_i M_n = 0$, or else there exists
an integer $j > i$ such that $A_j M_n \subsetneqq A_i M_n$.

For that purpose, we will follow a diagonal procedure involving `zigzags of
increasing amplitude':  If $M_1 = 0$, let $P_1$ be any left primitive ideal;
otherwise, pick $P_1$ left primitive with $P_1 M_1 \subsetneqq M_1$, and set $A_1 =
P_1$.  If
$A_1 M_2 = 0$, set $A_2 = A_1$; otherwise, choose a left primitive ideal $P_2$ with
$P_2 A_1 M_2 \subsetneqq A_1 M_2$, and set $A_2 = P_2 A_1$.  If $A_2 M_1 =
0$, set $A_3 = A_2$; otherwise, pick a left primitive ideal $P_3$ with
$P_3 A_2 M_1 \subsetneqq A_2 M_1$, and set $A_3 = P_3 A_2$.  In the next step,
we define $A_4 = P_4 A_3$, where $P_4$ is a left primitive ideal such that $P_4
A_3 M_2 \subsetneqq A_3 M_2$, unless $A_3 M_2 = 0$, and in the latter case we set
$A_4 = A_3$.  Now we move on to $M_3$, cutting down the size of the left module $A_4
M_3$ if the latter is nonzero, then we return to
$M_1$, $M_2$, $M_3$, from whence we move up to $M_4$,
following the pattern

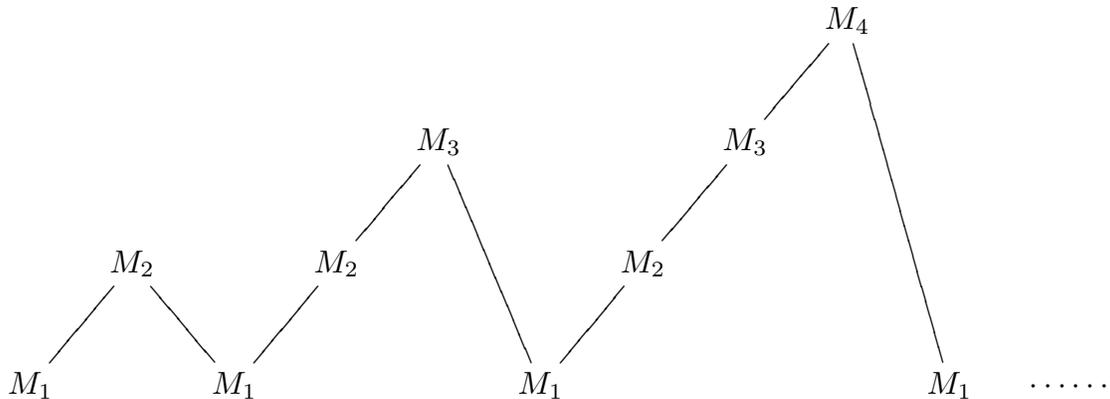
\ignore{
$$\xymatrixcolsep{1.0pc}
\xymatrix{
 &&&& &&&&M_4 \edge[dl]\edge[dddr]\\
 &&&&M_3 \edge[dl]\edge[ddr] &&&M_3 \edge[dl]\\
 &M_2 \edge[dl]\edge[dr] &&M_2\edge[dl] &&&M_2 \edge[dl]\\
M_1 &&M_1 &&&M_1 &&&&M_1 &\cdots\cdots
}$$
}

\noindent Inductively, this process yields a chain $A_1 \supseteq A_2 \supseteq \dots$
of ideals conforming to the above requirements.

Since all of the $A_n$ are finitely generated right ideals of $R$, the following
assignments yield
$p$-functors in the sense of \cite{\Zim}, i.e., subfunctors of the forgetful functor
$\Rmod
\rightarrow
\Zmod$ which commute with direct products:
$$F_n : \Rmod \rightarrow \Zmod, \quad \quad \quad X \mapsto
A_nX.$$
In the present situation, we are actually dealing with subfunctors $F_n$ of the
identity functor $\Rmod \rightarrow \Rmod$ taking finitely generated left $R$-modules
to finitely generated left $R$-modules.  Set
$M = \prod_{n \in
\NN} M_n$, and assume that, to the contrary of our claim, this product splits into
finitely generated direct summands, say 
$$\quad \quad M = \bigoplus_{l \in L} B_l, \tag\dagger$$
where each $B_l$ is a finitely generated left $R$-module.  An upgraded
version of Chase's Lemma \cite{\Huibiel, Lemma 11} (cf\. \cite{\Chase, Theorem
1.2} for the original result) then yields a natural number $n_0$ and a finite subset
$L_0
\subseteq L$ such that 
$$\quad \quad \prod_{n \ge n_0} F_{n_0}M_n
\quad \subseteq \quad \bigoplus_{l \in L_0} F_{n_0} B_l +
\bigcap_{k \in \NN}F_k M. \tag\ddagger$$
\medskip

{\it Case 1:}  $F_{n_0} M_n = 0$ for all but finitely many $n
\in \NN$

Set $N_0 = \{n \in \NN \mid F_{n_0} M_n = 0\}$ and $N_1  =  \NN
\setminus N_0$.  Factoring $F_{n_0} M$ out of both sides of
equality $(\dagger)$ yields
$$\prod_{n \in N_0} M_n \oplus \prod_{n \in N_1} \bigl(M_n /
F_{n_0} M_n \bigr) \quad  = \quad  \bigoplus_{l \in L}\bigl(B_l
/ F_{n_0} B_l \bigr).$$  
Observe that the summand $D:= \prod_{n \in N_1} \bigl(M_n /
F_{n_0} M_n \bigr)$ on the left-hand side of this last equality is a 
finitely generated $R$-module, since
$N_1$ is finite.  It is therefore contained in some finite subsum
$B: = \bigoplus_{l \in L_1} \bigl(B_l / F_{n_0} B_l \bigr)$ of the
right-hand direct sum, which is in turn finitely
generated over $R$.  In summary, we thus obtain
$$\prod_{n \in N_0} M_n \cong \bigoplus_{\l \in L \setminus
L_1} \bigl(B_l / F_{n_0} B_l \bigr) \oplus B / D,$$
where all summands on the right are finitely generated $R$-modules. 
Observe that all of the modules involved in this isomorphism
are annihilated by $A_{n_0}$, and so are modules over $R/ A_{n_0}$, the latter being an
Artin algebra by hypothesis.  But this is incompatible with the fact that our claim is
known to hold for such algebras:  Indeed, the infinite family $(M_n)_{n \in N_0}$
consists of pairwise non-isomorphic indecomposable
$R / A_{n_0}$-modules, since their $R/A_{n_0}$-structure coincides with their
$R$-structure.  This rules out the first case.
\medskip

{\it Case 2:}  $F_{n_0} M_n \ne 0$ for infinitely many $n \in
\NN$.

We denote by $N_1$ the infinite set $\{ n \in \NN
\mid n \ge n_0 \ \text{and}\ F_{n_0} M_n \ne 0 \}$.  Next we
factor the term $\bigcap_{k \in \NN} F_k M = \prod_{n \in \NN}\bigl( \bigcap_{k \in
\NN} F_k M_n \bigr)$ out of both sides of the inclusion labeled $(\ddagger)$ to
obtain
$$\prod_{n \ge n_0} \biggl( F_{n_0} M_n \bigm/\bigcap_{k \in
\NN} F_k M_n \biggr)  \subseteq \bigoplus_{l \in L_0}
\biggl( F_{n_0} B_l \bigm/  \bigcap_{k \in \NN}
F_k B_l \biggr).$$
In justifying this inclusion, keep in mind that $p$-functors automatically commute
with direct sums.  We note that the direct product on the left-hand side equals 
$$\prod_{n \in N_1} \biggl( F_{n_0} M_n \bigm/
\bigcap_{k \in \NN} F_k M_n \biggr).$$
Now $\bigcap_{k \in \NN} F_k M_n \subsetneqq F_{n_0} M_n$ for
all $n \in N_1$ by construction of the functors $F_k$, and consequently, the
left-hand side of the inclusion we just derived is not
finitely generated; this is due
to the fact that it contains the infinite direct sum $\bigoplus_{n \in N_1}
\biggl(F_{n_0} M_n \bigm/
\bigcap_{k \in \NN} F_k M_n \biggr)$. 
But the right-hand side of our inclusion {\it is} finitely generated, and we
have again reached an absurdity.

This completes the proof of the lemma. \qed \enddemo

\enddocument